\documentclass[12pt,letterpaper]{amsart}
\usepackage{euler, epic,eepic,latexsym, amssymb, amscd, amsfonts, xypic, url, color, epsfig, stmaryrd}

\newcommand{\tpoint}[1]{\vspace{3mm}\par \noindent \refstepcounter{subsection}{\bf \thesubsection.} 
  {\em #1. ---} }
\newcommand{\epoint}[1]{\vspace{3mm}\par \noindent \refstepcounter{subsection}{\thesubsection.} 
  {\em #1.} }
\newcommand{\bpoint}[1]{\vspace{3mm}\par \noindent \refstepcounter{subsection}{\bf \thesubsection.} 
  {\bf #1.} }

\newcommand{\kc}{k^s}

\newcommand{\Z}{\mathbb{Z}}

\newcommand{\Q}{\mathbb{Q}}

\newcommand{\G}{\mathbb{G}}

\newcommand{\Jac}{\operatorname{Jac}}

\newcommand{\Gal}{\operatorname{Gal}}
\newcommand{\Galk}{\operatorname{Gal}(k^s/k)}

\newcommand{\et}{\text{et}}

\newcommand{\Pic}{\operatorname{Jac}}

\newcommand{\N}{V}
\newcommand{\h}{f}
\newcommand{\Spec}{\operatorname{Spec}}

\newcommand{\hidden}[1]{\footnote{Hidden:  #1}}
\renewcommand{\hidden}[1]{}

\title[Massey products $\langle y,x,x,\ldots,x,x,y\rangle$]
{Massey products $\langle y,x,x,\ldots,x,x,y\rangle$ in Galois cohomology via rational points}
\author[Wickelgren]{Kirsten Wickelgren}\thanks{The author is partially supported by NSF DMS-1406380.}

\address{School of Mathematics, Georgia Institute of Technology}

\email{wickelgren@post.harvard.edu}

\date{January 17, 2016}

\subjclass[2010]{Primary 55S30, Secondary 11S25, 14H30. }

\keywords{Massey products, Galois cohomology, Section conjecture, Lower central series, Group cohomology}

\begin{document}

\begin{abstract}
For $x$ an element of a field other than $0$ or $1$, we compute the order $n$ Massey products $$\langle (1-x)^{-1}, x^{-1}, \ldots, x^{-1},  (1-x)^{-1} \rangle$$ of $n-2$ factors of $x^{-1}$ and two factors of $(1-x)^{-1}$ by embedding $\mathbb{P}^1 - \{0,1,\infty\}$ into its Picard variety and constructing $\Gal(\kc/k)$ equivariant maps from $\pi_1^{\et}$ applied to this embedding to unipotent matrix groups. This method produces obstructions to $\pi_1$-sections of $\mathbb{P}^1 - \{0,1,\infty\}$, partial computations of obstructions of Jordan Ellenberg, and also computes the Massey products $$\langle x^{-1} , (-x)^{-1}, \ldots, (-x)^{-1}, x^{-1} \rangle.$$ 

\end{abstract}

\maketitle

\section{Introduction} 

A smooth curve $X$ over a field $k$ with a choice of basepoint can be embedded into its generalized Jacobian, $\Pic X$. Given a rational point $\tilde{p}$ of $X$, there is a commutative diagram $$\xymatrix{ & X \ar[d] \\ \Spec k \ar[r] \ar[ur]^{\tilde{p}} & \Pic X.}$$ Applying the \'etale fundamental group functor $\pi_1^{\et}$, with the same chosen basepoint, to this diagram, we obtain the diagram $$\xymatrix{ &\pi_1^{\et} X \ar[d] \\ \pi_1^{\et} \Spec k \ar[r] \ar[ur] & \pi_1^{\et}\Pic X.}$$ Given a rational point $p$ of $\Pic X$, there is an associated map  \begin{equation}\label{pi1SpecktoPic} \pi_1^{\et} \Spec k \to   \pi_1^{\et}\Pic X.\end{equation} The absence of a lift of \eqref{pi1SpecktoPic} through  the map \begin{equation}\label{pi1XtoPic} \pi_1^{\et} X \to   \pi_1^{\et}\Pic X\end{equation} is an obstruction to the existence of a rational point of $X$ mapping to $p$. 

For $X =\mathbb{P}^1_k - \{0,1,\infty\}$, we construct obstructions to lifting a map \eqref{pi1SpecktoPic} through the map \eqref{pi1XtoPic} and compute them to be higher order Massey products. Since rational points themselves are unobstructed, this produces computations of certain Massey products. 

More specifically, these computations of Massey products are as follows. Choose an integer $n \geq 3$ and let $\Sigma$ denote a set of primes not dividing $n!$ and not including the characteristic of $p$. For a profinite group $\pi$, let $\pi^{\Sigma}$ denote the maximal pro-$\Sigma$ quotient of $\pi$. Let $\chi: \Galk \to (\Z^{\Sigma})^*$ denote the cyclotomic character (see \ref{Applications_section_notation_subsection}), and let $\Z^{\Sigma}(\chi)$ denote the module over $\Galk$ whose underlying abelian group is $\Z^{\Sigma}$ and where the $\Galk$ action is via $\chi$. Identify elements $x$ of $k^*$ with their images under the Kummer map $$k^* \to H^1(\Galk, \Z^{\Sigma}(\chi)),$$ recalled below in \ref{Kummer_map_recall}.  The definition of the $n$th order Massey product of $n$ elements of $H^1(\Galk, \Z^{\Sigma}(\chi))$ is recalled in \ref{Massey_prod_subsection}. There is a cocycle $\h: \Galk \to \Z^{\Sigma}(\chi^{n-1})$ associated to the $\Galk$-action on $\pi_1^{\et, \Sigma}(\mathbb{P}^1_{\kc} - \{0,1,\infty\})$; it is described in Theorem \ref{Main_thm_obstruction_section}. Here, the fundamental group is taken with respect to the rational tangential base point $0 + \epsilon$. Tangential base points and the notation used here are recalled in \ref{rational_tangential_basepoints_subsection} and \ref{notation_for_tangential_base_points}. 

\tpoint{Theorem}\label{Intro_cor_Cor_MPGA}{\em Let $x$ be an element of $k^* -\{1\}$. The $n$th order Massey product $$\langle (1-x)^{-1}, x^{-1}, \ldots, x^{-1}, (1-x)^{-1} \rangle $$ contains $$0 \text { and }-[\h] \cup (1-x)^{-1}.$$ Let $x$ be an element of $k^*$. The $n$th order Massey product $$\langle -x^{-1}, x^{-1}, \ldots, x^{-1}, -x^{-1} \rangle $$ contains $$0 \text { and }-([\h] \cup (-x)^{-1}).$$}

In Theorem \ref{Intro_cor_Cor_MPGA}, the Massey product has $(n-2)$ factors of $x^{-1}$ and two factors of $(1-x)^{-1}$. Theorem \ref{Intro_cor_Cor_MPGA} is proved as Corollary \ref{Cor_MPGA} below. Recall  that associated to every element of an $n$th order Massey product there is a defining system (see \ref{Massey_prod_subsection}). The defining systems implicit in Theorem \ref{Intro_cor_Cor_MPGA} are specified in Corollary \ref{Cor_MPGA} in terms of an $(n-1)$-nilpotent quotient of the map $ \pi_1^{\et} \Spec k \to \pi_1^{\et}(\mathbb{P}^1_k - \{0,1,\infty\})$ associated to the rational point or tangential point corresponding to $x$ in $k^* -\{1\}$ or $k^*$.

The obstructions giving rise to Theorem \ref{Intro_cor_Cor_MPGA} are as follows. Let $\pi =  \pi_1^{\et, \Sigma}(\mathbb{P}^1_{\kc} - \{0,1,\infty\}, 0 + \epsilon)$ and $G = \Galk \cong \pi_1^{\et} \Spec k$. The profinite group $\pi$ has a continuous action of $G$. Let $$\pi=[\pi]_1 \supseteq [\pi]_2 \supseteq [\pi]_3 \ldots$$ denote its lower central series (see \ref{notation_subsection}). There is a canonical isomorphism $\pi_1^{\et, \Sigma}(\Pic \mathbb{P}^1_{\kc} - \{0,1,\infty\}) \cong \pi/[\pi]_2$ such that applying $\pi_1^{\et, \Sigma}(- \otimes \kc)$ to the embedding $X \to \Pic X$ gives the abelianization map. This is recalled in \ref{subsection_pi1_sectionsP1minus}. For a quasi-compact, geometrically connected scheme $X$ over $k$ with a rational basepoint (which can be tangential), the homotopy exact sequence \cite[IX Th. 6.1, Co. 4.11]{sga1} gives a canonical isomorphism $\pi_1(X) \cong \pi_1(X_{\kc}) \rtimes G$.  Choose $n \geq 3$. Suppose given a map \eqref{pi1SpecktoPic}, which then yields a homomorphism over $G$\begin{equation*} G \to \pi/[\pi]_2 \rtimes G \end{equation*} or equivalently a cocycle \begin{equation}\label{Gtopiab_map} G \to \pi/[\pi]_2.\end{equation}

$\pi$ is isomorphic to the pro-$\Sigma$ completion of the free group on two generators $x$ and $y$,  and the $G$-action on $\pi$ is of the form $$g(x) =x^{\chi(g)} $$ $$g(y) =\mathfrak{f}(g)^{-1}y^{\chi(g)}\mathfrak{f}(g) $$ where $\mathfrak{f}: G \rightarrow [\pi]_2$ is a cocycle taking values in the commutator subgroup of $\pi$. This is recalled in \ref{description_pi_and_G_action}. Let $U_{n+1}$ denote the upper triangular $(n+1) \times (n+1)$ square matrices with entries in $\Z^{\Sigma}$, diagonal entries $1$, and let $\overline{U}_{n+1}$ denote the quotient of $U_{n+1}$ by those matrices whose only non-zero entry is the one in the upper right corner. There is a $G$-action on $U_{n+1}$ and $\overline{U}_{n+1}$ such that weights of commutators are compatible with those of the fundamental group of $\mathbb{P}^1 - \{0,1,\infty\}$. This action is defined explicitly in \ref{notation_subsection}. The map taking a matrix to the entries just above the diagonal defines a quotient map $\overline{U}_{n+1} \to \oplus_{k=1}^n \Z^{\Sigma}(\chi)$. Note that $\pi/[\pi]_2$ is isomorphic to $\pi/[\pi]_2 \cong \Z^{\Sigma}(\chi)^2$. We show that the $G$-equivariant map \begin{equation}\label{piab_to_Zsigmachin_map}\pi/[\pi]_2 \to  \oplus_{k=1}^n \Z^{\Sigma}(\chi)\end{equation} defined $$y \mapsto 1 \oplus 0 \oplus 0 \ldots 0 \oplus 0 \oplus 1$$  $$x \mapsto 0 \oplus 1 \oplus 1 \ldots \oplus 1 \oplus 1 \oplus 0 $$ lifts to a $G$-equivariant map \begin{equation}\label{pitoUn+1_into}\pi \to \overline{U}_{n+1}.\end{equation} By composing with the map \eqref{piab_to_Zsigmachin_map}, the map \eqref{Gtopiab_map} yields a cocycle $$G \to \oplus_{k=1}^n \Z^{\Sigma}(\chi).$$ If there is a lift of \eqref{pi1SpecktoPic} through \eqref{pi1XtoPic}, then composing the corresponding cocycle with \eqref{pitoUn+1_into} yields a cocycle $G \to \overline{U}_{n+1}$. Such a map is the same data as a defining system for an $n$th order Massey product. We evaluate the corresponding Massey product in Theorem \ref{Main_thm_obstruction_section}, which is the main theorem of the paper. Thus, if there is a lift of \eqref{pi1SpecktoPic} through \eqref{pi1XtoPic}, then an $n$th order Massey product must be as computed, and this is the obstruction to be elaborated in this paper. 

The specification of the defining system allows the obstruction to be combined with others. For example, the obstruction here for $n=3$ combined with an obstruction in \cite{Kyoto} were studied in \cite{PIA}, where they are used to find maps \eqref{pi1SpecktoPic} which do not lift through \eqref{pi1XtoPic}.

The obstruction given here factors through obstructions of Jordan Ellenberg \cite{Ellenberg_2_nil_quot_pi} which we recall in \ref{subsection_Ellenberg_obstructions}. It gives a partial computation of Ellenberg's obstructions in terms of Massey products. See Remark \ref{main_theorem_as_partial_computation_delta_n} and \eqref{mu_*delta_n_comp}. 

Motivation for studying such obstructions includes Grothendieck's section conjecture which predicts that for a smooth curve $X$ of negative Euler characteristic over a number field $k$, the rational points and rational tangential points at infinity (see \ref{rational_tangential_basepoints_subsection} for the definition of a rational tangential point at infinity) are in natural bijection with the conjugacy classes of sections of the homotopy exact sequence $$1 \to \pi_1^{\et}(X_{\overline{k}}) \to \pi_1^{\et}(X) \to \Gal(\overline{k}/k) \to 1,$$ where the conjugacy class of a section $\tilde{p}: \Gal(\overline{k}/k) \to \pi_1(X)$ is the set of sections of the form $g \mapsto \gamma \tilde{p}(g) \gamma^{-1}$ with $\gamma \in \pi_1^{\et}(X_{\overline{k}})$. Call the conjugacy class of a section such as $\tilde{p}$ a {\em $\pi_1$-section}. 

The curve $\mathbb{P}^1 - \{0,1,\infty\}$ is an important example for studying the section conjecture by Belyi's theorem, which says that every smooth proper curve over a number field is a branched cover of $\mathbb{P}^1$ only branched over $\{0,1,\infty\}$ \cite{Belyi}, combined with the fact that the section conjecture behaves well with respect to covering space maps \cite[Ch 9 Propositions 110 and 104]{Stix_evidence}. To study the section conjecture for $\mathbb{P}^1 - \{0,1,\infty\}$, we compare the union of the rational points and rational tangential points at infinity with the $\pi_1^{\et}$-sections. As the former are completely understood, it is the latter that need to be computed to equal the former. 

The $\pi_1$-sections of $\Pic (\mathbb{P}^1 - \{0,1,\infty\})$ in characteristic $0$ can be computed to be the set of elements of the form $[p] = [p]_x \oplus [p]_y$ with $[p]_x$ and $[p]_y$ in $H^1(\Galk, \Z^{\wedge}(\chi))$, where $\Z^{\wedge}$ denotes the profinite completion of $\Z$. This is discussed further in \ref{subsection_pi1_sectionsP1minus}. In \cite{PIA} and \cite{Kyoto}, the problem of determining which $\pi_1$-sections of $\Pic (\mathbb{P}^1 - \{0,1,\infty\})$ lift to $\mathbb{P}^1 - \{0,1,\infty\}$ was studied. We continue this study here: as a corollary of the main theorem of this paper, we have:

\tpoint{Theorem}\label{Intro_version_Cor_pi1-sections_Pic_lift_obstruction}{\em A $\pi_1^{\et,\Sigma}$-section $p$ of $\Pic (\mathbb P_k^1 - \{0,1,\infty \})$ which lifts to a $\pi_1^{\et,\Sigma}$-section of $\mathbb P_k^1 - \{0,1,\infty \}$ satisfies the condition that the Massey product $$\langle-[p]_y, -[p]_x, \ldots, -[p]_x, -[p]_y  \rangle$$ contains  $$0 \text { and }-[\h] \cup [p]_y.$$ }

This is proven as Corollary \ref{Cor_pi1-sections_Pic_lift_obstruction} below. Because the defining system is specified in Theorem \ref{Main_thm_obstruction_section}, we have a similar specification in Corollary \ref{Cor_pi1-sections_Pic_lift_obstruction} as well. 

\epoint{Remark} The position of the two occurrences of $-[p]_y$ in the Massey product of Theorem \ref{Intro_version_Cor_pi1-sections_Pic_lift_obstruction} must be as written for the method of this paper to give a computation of the Massey product. The reason for this is that the proof that the lift \eqref{pitoUn+1_into} of\eqref{piab_to_Zsigmachin_map} is equivariant relies on the given position. For example, the given position is essential in Lemma \ref{N_normal_B_in_center}.

\epoint{Remark} In \ref{Intro_cor_Cor_MPGA},  \ref{Intro_version_Cor_pi1-sections_Pic_lift_obstruction} and the main theorem, it is possible to pro-prime-to-$p$ complete and invert all primes less than or equal to $n$, instead of using pro-$\Sigma$ completions. See Remarks \ref{main_thm_mod_remark} and \ref{corMPGA_mod_remark}.

\subsection{Relation with other work} Guillou \cite{Guillou} and Sharifi \cite{Sharifi} have interesting overlapping and related results. Guillou shows that all Massey products of order $\leq 4$ in motivic cohomology of the elements $x$ and $1-x$ are defined and contain $0$, except possibly for $\langle x,x,1-x,1-x\rangle$ and $\langle x, 1-x, x, 1-x \rangle$. The \'etale realization of this result should produce the analogous computation in Galois cohomology. Sharifi shows that many Massey products in Galois cohomology of the form $\langle x,x, \ldots, x,y \rangle$ are defined and contain $0$ under certain hypotheses. These are both consistent with the results here and in \cite{PIA}, \cite{Kyoto}.

{\em Acknowledgments:} I thank the referee for useful comments. 

\section{$n$-Nilpotent obstruction}

The goal of this section is to prove Theorem \ref{Main_thm_obstruction_section}, computing a type of Massey product. The material until \ref{subsection_some_defining_systems} is expository or serves to establish notation.

\subsection{Notation}\label{notation_subsection}
In this section, $\Sigma$ is a set of primes of $\Z$, and $\Z^{\Sigma}$ is the pro-$\Sigma$ completion of $\Z$.

$\pi$ is the pro-$\Sigma$ completion of the free group $\langle x, y \rangle$ on two generators $x$ and $y$.  $G$ is a group or a profinite group which acts continuously on $\pi$ by the formula \begin{align}\label{g_action_pi} &g(x)  =x^{\chi(g)} \nonumber \\ g(y) & =\mathfrak{f}(g)^{-1}y^{\chi(g)}\mathfrak{f}(g),\end{align} where $\chi: G \to (\Z^{\Sigma})^*$ is a continuous homomorphism. 

For any integer $m$, let $\Z^{\Sigma}(\chi^m)$ denote the group $\Z^{\Sigma}$ with the continuous action of $G$ where $g$ in $G$ acts by multiplication by $\chi(g)^m$.

For elements $\gamma_1$,$\gamma_2$ of a group, the commutator $[\gamma_1, \gamma_2]$ is defined to be $[\gamma_1, \gamma_2] = \gamma_1 \gamma_2  \gamma_1^{-1} \gamma_2^{-1}$. Given a group $\Gamma$, the {\em lower central series} of $\Gamma$ is denoted $$\Gamma = [\Gamma]_1 \supseteq [\Gamma]_2 \supseteq [\Gamma]_2 \supseteq \ldots,$$ and is defined by setting $[\Gamma]_{n+1} = \overline{[\Gamma, [\Gamma]_n]}$ to be the closure of the subgroup generated by commutators of an element of $
\Gamma$ with an element of $[\Gamma]_n$. 

For a profinite abelian group $A$ with a continuous action of $G$, let $C^n(G,A)$ denote the group of in homogenous $n$-cochains for $n=1,2,\ldots$ as in \cite[p. 14]{coh_num_fields}. Let $D: C^n(G,A) \to C^{n+1}(G,A)$ denote the differential, and $H^n(G,A)$ denote the corresponding cohomology groups. Suppose $A'$ is also a profinite abelian group with a continuous action of $G$. For $a \in C^1(G, A)$ and $a' \in C^1(G, A')$, let $a \cup a'$ denote the cup product $a \cup a' \in C^{2}(G, A \otimes A')$  $$(a \cup a')(g_1,g_2) = a(g_1) \otimes (g_1a'(g_{2})).$$ This product induces a well defined map on cohomology. 

For $A$ a profinite group with a continuous action of $G$ which is not necessarily assumed to be abelian, let $C^1(G,A)$ denote the set of functions $a:G \to A$ such that $a(gh) = a(g) (g a(h))$, and let $H^1(G,A)$ denote the pointed set of equivalence classes of element of $C^1(G,A)$ where $a,a'$ are equivalent if and only if there exists $\gamma$ in $A$ such that $a'(g) = \gamma^{-1} a(g) g\gamma$ for all $g$ in $G$ \cite[p. 16]{coh_num_fields}.

Let $a_{ij}$ be the function taking a matrix to its $(i,j)$-entry.

As usual, a matrix is said to be {\em upper triangular} if all the entries below the main diagonal are $0$, or in other words, if $a_{i,j} = 0$ for $j < i$. A matrix is said to be {\em strictly upper triangular matrix} if the matrix is upper triangular and the diagonal entries are $0$, or in other words if $a_{i,j} = 0$ for $j \leq i$.

For a matrix $M$, let $M_{ij}$ denote the matrix entry in the $i$th row and $j$th column. Let $E_{ij}$ be the matrix such that $a_{ij}(E_{ij}) = 1$ and for all $(k,l) \neq (i,j)$, we have $a_{kl}( E_{ij}) = 0$. Let $1$ denote the identity matrix.

Let $U_{n+1}$ denote the multiplicative group of $(n+1) \times (n+1)$ upper triangular matrices with coefficients in $\Z^{\Sigma}$ whose diagonal entries are $1$. ``U" stands for unipotent, not unitary.  Give $U_{n+1}$ a $G$-action by defining $g M$ via $$a_{ij} (g M )= \chi(g)^{j-i} a_{ij}(M),$$ for all $M \in U_{n+1}$ and $g$ in $G$. Note that $U_{n+1}$ is a pro-$\Sigma$ group; it is isomorphic to the inverse limit of matrix groups with coefficients in various cyclic abelian groups whose orders are only divisible by primes in $\Sigma$. Given $M$ in $U_{n+1}$ and $N$ in $\Z^{\Sigma}$, there is a natural definition of $M^N$; $M$ can be expressed as an inverse limit of matrices of finite order $M= \varprojlim_j M_j$ and $M^N = \varprojlim_j M_j^N$, where $M_j^N$ is well-defined because a power of $M_j$ only depends on the exponent modulo the order of $M_j$. The $G$-action on $U_{n+1}$ is continuous. 

There is a $G$-equivariant inclusion $\Z^{\Sigma}(\chi^n) \rightarrow U_{n+1}$ sending $m$ in $\Z^{\Sigma}$ to the matrix $1+ mE_{1,n+1}$.  This inclusion gives rise to central extensions \begin{align}
&1 \to \Z^{\Sigma}(\chi^n) \to U_{n+1}  \to \overline{U}_{n+1} \to 1 \nonumber \\
\label{Z_Un+1_barUn+1}
1 \to &\Z^{\Sigma}(\chi^n) \to U_{n+1}  \times G \to \overline{U}_{n+1} \rtimes G \to 1.
\end{align} \noindent where $\overline{U}_{n+1}$ is defined as the quotient $U_{n+1}/\Z^{\Sigma}(\chi^n)$. Note that $\overline{U}_{n+1}$ is also a profinite group.

\subsection{Massey Products}\label{Massey_prod_subsection}

We recall the definition and some well-known properties of Massey products. Let $z_1, \ldots z_n$ be elements of $H^1(G, \Z^{\Sigma}(\chi) )$. A {\em defining system} for the $n$th order Massey product $\langle z_1, z_2, \ldots, z_n \rangle$ consists of cochains $Z_{i,j}$ in $C^1( G, \Z^{\Sigma}(\chi^{j-i}))$ for $i$ and $j$ in $\{1,2,\ldots, n+1 \}$ such that $i<j$ and $(i,j) \neq (1,n+1)$ satisfying 

\begin{itemize}
\item $Z_{i,i+1}$ represents $z_i$,
\item $D Z_{i,j} = \sum_{r=i+1}^{j-1} Z_{i,r} \cup Z_{r,j}$ for $i+1<j$.
\end{itemize} 

The $n$th order Massey product $\langle z_1, z_2, \ldots, z_n \rangle$ is defined when a defining system exists, and the {\em Massey product relative to this defining system} is the cohomology class in $H^2(G, \Z^{\Sigma}(\chi^n) )$ determined by the cocycle $$ \sum_{r=2}^{n} Z_{1,r} \cup Z_{r,n+1}.$$ The {\em Massey product} $\langle z_1,\ldots z_n \rangle $ is the set of all cohomology classes associated to such defining systems.

\epoint{Example} See \cite[\S 2]{Dwyer} or \cite[2.4]{Kyoto}. The functions $a_{i,j}$ determine cochains in $C^1(\overline{U}_{n+1} \rtimes G, \Z^{\Sigma}(\chi^{j-i}))$. It follows from the definition of matrix multiplication and the semi-direct product that the boundary of $a_{i,j}$  is computed by \begin{equation}\label{Dai,j}D (a_{i,j}) = - \sum_{r=i+1}^{j-1} a_{i,r} \cup a_{r,j}.\end{equation} Thus the cochains $-a_{i,j}$ for $i$ and $j$ in $\{1,2,\ldots, n+1 \}$ with $i<j$ and $(i,j) \neq (1,n+1)$ form a defining system for the $n$th order Massey product $$\langle -a_{1,2}, -a_{2,3}, \ldots, -a_{n,n+1}\rangle.$$ It is convenient to introduce the following notation for this particular example. 

\bpoint{Definition}\label{def_defining_system_determined_by_theta} Let $\theta$ be a homomorphism $\theta: G \to \overline{U}_{n+1} \rtimes G$ such that the composition of $\theta$ with the quotient map $\overline{U}_{n+1} \rtimes G \to G$ is the identity. Then the cochains \begin{equation}\label{defining_system_determined_theta} -a_{i,j} \theta \text{ for } i,j \in \{1,2,\ldots, n+1 \},~~ i<j \text{ and } (i,j) \neq (1,n+1) \end{equation} form a defining system for the $n$th order Massey product $$\langle -a_{1,2} \theta, -a_{2,3} \theta, \ldots, -a_{n,n+1} \theta\rangle.$$ We will call the defining system \eqref{defining_system_determined_theta} the {\em defining system determined by $\theta$.} Up to sign conventions, and adding a non-trivial action, this definition is implicit in \cite[Theorem 2.6]{Dwyer}. 

\subsection{Ellenberg's obstructions}\label{subsection_Ellenberg_obstructions} We recall the definition of Jordan Ellenberg's obstructions. In the context of Grothendieck's section conjecture and studying the rational points of curves, Ellenberg suggested successively studying the images of $$H^1(G, \pi/[\pi]_n) \to H^1(G, \pi/[\pi]_2) $$ using the boundary maps $$\delta_n:  H^1( G, \pi/[\pi]_n) \to H^2( G, [\pi]_n/[\pi]_{n+1})$$ coming from the central extensions $$1 \to [\pi]_n/[\pi]_{n+1} \to \pi/ [\pi]_{n+1} \to  \pi/ [\pi]_n \to 1.$$

\subsection{Some defining systems.}\label{subsection_some_defining_systems} We now construct the defining systems associated to our obstructions.

\bpoint{Assumption}\label{assumption_n>=3}{\em Fix an integer $n$ with  $n \geq 3$.} 

Let $A$ in $U_{n+1}$ denote the matrix such that $$a_{i,j} A =  \begin{cases}
 \frac{1}{(j-i)!} &1<i<j<n+1 \\
 1 & $j=i$ \\
 0 & \mbox{otherwise}
 \end{cases}.$$ For example, when $n=4$ we have $$A = \left( \begin{array}{ccccc}
1 & 0 & 0 & 0 & 0 \\
0 & 1 & 1 & \frac{1}{2} & 0 \\
0 & 0 & 1 & 1 & 0 \\ 
0 & 0 & 0 & 1 & 0 \\
0& 0 & 0 & 0 & 1 \end{array} \right)$$ By a slight abuse of notation, we will also let $A$ denote its image in $\overline{U}_{n+1}.$

Let $\N \subset \overline{U}_{n+1}$ denote the set of those elements $\overline{M}$ such that $a_{ij} (\overline{M}) = 0$ for $i>1$ and $j< n+1$ and $i \neq j$ \begin{equation} \N = \{\overline{M} \in \overline{U}_{n+1} : a_{ij} (\overline{M}) = 0 \text{ if } i>1, j< n+1,\text{ and } i \neq j \} .\end{equation}

Let $B = 1 + E_{1,2}+E_{n,n+1}$. For example, when $n=4$ we have $$B = \left( \begin{array}{ccccc}
1 & 1 & 0 & 0 & 0 \\
0 & 1 & 0 & 0 & 0 \\
0 & 0 & 1 & 0 & 0 \\ 
0 & 0 & 0 & 1 & 1 \\
0& 0 & 0 & 0 & 1 \end{array} \right)$$ By a slight abuse of notation, we will also let $B$ denote its image in $\N$. 

In the notation of \cite[ p. 586 (7)]{Kyoto}, $A = \varphi_J(\gamma_1)$ and $B= \varphi_J(\gamma_2)$ for $J: \{ 1,2, \ldots, n\} \to \{1,2\}$ the function such that $J(1)= J(n) = 2 $ and $J(2) = J(3) = \ldots = J(n-1) = 1$.

\bpoint{Definition of $\varphi$}
{\em Let $\varphi: \pi \to \overline{U}_{n+1}$ be the continuous homomorphism defined by $\varphi (x) = A$ and $\varphi(y) = B$.} 

We will show that $\varphi$ is a $G$-equivariant homomorphism using the following lemmas.

\tpoint{Lemma}\label{AandB_raised_to_N}{\em For all positive integers $N$, we have \begin{itemize}
\item $a_{i,i+j} (A^N) = N^j a_{i, i+j} A$.
\item $a_{i,i+j} (B^N) = N^j a_{i, i+j} B$.
\end{itemize} }

\begin{proof}
In \cite[Lemma 1]{Kyoto}, it was shown that  $a_{i,i+j} (A_{\ell}^N) = N^j a_{i, i+j} A_{\ell}$, where $A_{\ell}$ is the matrix in $U_{\ell + 1}$ defined by $a_{i, i+j} A_{\ell} = \frac{1}{j!}$ for $j>0$. Since $A$ and $B$ are block diagonal matrices whose blocks are of the form $A_{\ell}$ for various $\ell$'s, the lemma follows.
\end{proof}

\tpoint{Lemma}\label{N_normal_B_in_center}{\em 
$\N$ is a normal subgroup which contains $B$ in the center.}

\begin{proof}
Say that an upper triangular $(n+1) \times (n+1)$ matrix $C$ satisfies Condition \eqref{condition_star} if 
\begin{equation}\label{condition_star}
C_{ij} = 0 \text{ for } i>1 \text{ and } j< n+1.
\end{equation}
Let $C$ and $D$ be strictly upper triangular matrices, with $C$ satisfying Condition \eqref{condition_star}. Then $$a_{i,j}(CD) = \Sigma_{k = 1}^{n+1} C_{ik}D_{kj} = \Sigma_{i<k<j} C_{ik}D_{kj},$$ where the second equality holds because $C$ and $D$ are upper triangular matrices such that $C_{ii} = D_{ii} = 0.$ Choose $i>1$ and $j<n+1$. Then for any $k$ such that $i<k<j$, we have $k < n+1$, whence $C_{ik} = 0$ by Condition \eqref{condition_star}, from which it follows that $a_{i,j} (C D) = 0$. 

Any element of $\N$ is the image of a matrix of the form $1 + C$ with $C$ a strictly upper triangular matrix satisfying Condition \eqref{condition_star}. Since $(1 + C)(1+D) = 1 + C + D + CD$, it follows that $\N$ is closed under multiplication. Since $(1+C)^{-1} =1 + \Sigma_{n=1}^{\infty} (-1)^n C^n $, it also follows that $\N$ is closed under taking inverses, from which it follows that $\N$ is a subgroup. 

Let $C$ and $D$ again be strictly upper triangular matrices with $C$ satisfying Condition \eqref{condition_star}. Note that $$(1+D) (1+ C) (1+ D)^{-1} = (1+ D)(1+C)(1 + \Sigma_{n=1}^{\infty} (-1)^n D^n ).$$ By the above, we have that $D C D^n$ and $C D^n$ satisfy Condition \eqref{condition_star}. It follows that the image of $(1+D) (1+ C) (1+ D)^{-1}$ in $\overline{U}_{n+1}$ lies in $\N$. Since any element of $\overline{U}_{n+1}$ is the image of an element of $U_{n+1}$ of the form $1+D$, it follows that $\N$ is a normal subgroup.

Let $C$ again be a strictly upper triangular matrix satisfying Condition \eqref{condition_star}. To show that $B$ is in the center of $\N$, it suffices to show that $a_{i,j}(B(1+C) - (1+ C)B) = 0 $ for $(i,j) \neq (1,n+1)$, and this is what we do. Since $B (1 + C) = (1 + E_{1,2}+E_{n,n+1}) (1 + C )$ and $(1+C)B$ can be expressed similarly, it follows that \begin{equation}\label{Bcommutator} B(1+C) - (1+ C)B  = E_{1,2} C - C E_{1,2} +  E_{n,n+1} C - C E_{n,n+1}.\end{equation}  Since $\N$ is a subgroup, we have that $B(1+C)$ and $(1+C) B$ are upper triangular matrices whose diagonal entries are $1$ and such that $a_{i,j}=0$ when $i \neq j$, $i>1$ and $j<n+1$ . Thus $$a_{i,j}(B(1+C) - (1+ C)B) = 0 $$ unless $i = 1$ or $j = n+1$. By \eqref{Bcommutator}, it follows that \begin{align*}a_{1,j}(B(1+C) - (1+ C)B) = a_{1,2}(E_{1,2}) C_{2,j} - a_{i,j}(C E_{n,n+1} ) \\ = \begin{cases} 0 &\mbox{if } j<n+1 \\ 
 C_{2,n+1} - C_{1,n}& \mbox{if } j=n+1 \end{cases},\end{align*} where the first equality follows because $C$ is a strictly upper triangular matrix, and the second follows because $C_{2,j} = 0$ for $j \neq n+1$. 
 
For $i >1$, note that $a_{i,n+1} (E_{1,2} C) = 0$, $a_{i,n+1} (C E_{1,2} ) = 0$. Because $C$ is a strictly upper triangular matrix, it also follows that $a_{i,n+1} (E_{n,n+1} C) = 0$.  Because $C$ satisfies Condition \eqref{condition_star} and $i>1$ and $n< n+1$, we have that $a_{i,n+1} (C E_{n,n+1}) = 0$. By \eqref{Bcommutator}, it follows that $a_{i,n+1}(B(1+C) - (1+ C)B) = 0$.

This completes the proof of the claim.
\end{proof}

\tpoint{Lemma}\label{varphi(commutator)_in_N}{\em For any $\gamma$ in $[\pi]_2$, we have $\varphi(\gamma)$ is in $\N$.}

\begin{proof}
$[\pi]_2$ is topologically generated by elements of the form $$[\cdots[[x,y],z_1],z_2, \ldots],z_k],$$ where $z_i$ is either $x$ or $y$ and $k=0,1,2,\ldots$. By Lemma \ref{N_normal_B_in_center}, $\N$ is a normal subgroup. By the definition of $\varphi$, we have that $\varphi(y)$ is contained in $\N$. It follows that $[\varphi(x),\varphi(y)]$ is contained in $\N$, and more generally that $[\cdots[[\varphi(x),\varphi(y)],\varphi(z_1)],\varphi(z_2), \ldots],\varphi(z_k)]$ is contained in $\N$, proving the lemma.

\end{proof}

To compute our obstructions as Massey products, we define the following lift of $\varphi$.

\bpoint{Definition of $\varphi'$}
{\em Define $\varphi': \pi/[\pi]_{n+1} \to U_{n+1}$ to be the continuous homomorphism defined by $\varphi'(x) = A$ and $\varphi'(y) = B$.}

Note that $\varphi'$ is not claimed to be an equivariant homomorphism. 

\tpoint{Lemma}\label{gvarphi'varphi'gdif}{\em $g \varphi'(x)= \varphi'(x^{\chi(g)}) = \varphi' (g x)$ and $g \varphi' (y) = \varphi'(y^{\chi(g)})$}

\begin{proof}
Lemma \ref{AandB_raised_to_N} implies that $a_{i,i+j}(A^{\chi(g)}) = \chi(g)^j a_{i,i+j}(A)$.  By the definition of the $G$-action on $U_{n+1}$, it follows that $$ A^{\chi(g)} = g A .$$ By definition of $\varphi'$, it follows that $ A^{\chi(g)} = g\varphi'(x)$. Since $\varphi'$ is a homomorphism, it follows that $\varphi'(x^{\chi(g)}) = g\varphi'(x)$.

The same reasoning shows that $\varphi'(y^{\chi(g)}) = g \varphi'(y)$.

By the definition of the $G$-action on $\pi$, $g x = x^{\chi(g)}$, from which is follows that $\varphi'(gx) =  \varphi' (x^{\chi(g)}) $. This completes the proof of the lemma.
\end{proof}

\tpoint{Proposition}{\em $\varphi: \pi \to \overline{U}_{n+1}$ is $G$-equivariant.}

\begin{proof}
By the Lemma \ref{gvarphi'varphi'gdif}, we have $\varphi (g x) = g \varphi(x)$ and $\varphi(y^{\chi(g)}) = g \varphi (y) $. By \eqref{g_action_pi}, we have that $$\varphi (g y) = \varphi( \mathfrak{f}(g))^{-1}\varphi (y)^{\chi(g)} \varphi (\mathfrak{f}(g)).$$ Since $\mathfrak{f}(g) \in [\pi]_2$,   we have that $ \varphi (\mathfrak{f}(g))$ is in $\N$ by Lemma \ref{varphi(commutator)_in_N}. By Lemma \ref{N_normal_B_in_center},  $B = \varphi (y)$ is in the center of $\N$. It follows that $B^{\chi(g)}$ is in the center of $\N$. Thus $$\varphi( \mathfrak{f}(g))^{-1}\varphi (y)^{\chi(g)} \varphi (\mathfrak{f}(g))=\varphi (y)^{\chi(g)} = \varphi(y^{\chi(g)}).$$ Thus $\varphi (g y) = g \varphi(y).$ This proves that $\varphi$ is $G$-equivariant as claimed. 
\end{proof}

Recall that $n$ is assumed to be greater than $2$.

\tpoint{Lemma}\label{[B,C]CinN}{\em Suppose $C$ in $U_{n+1}$ maps to an element of $\N$ under the quotient map $U_{n+1} \to \overline{U}_{n+1}$. Then $[B, C]$ is in the kernel of this quotient map and $a_{1,n+1} [B, C] = a_{2,n+1}(C) - a_{1,n} C$.}

\begin{proof}
$C$ may be written as $C = 1 +D$, where $D$ is strictly upper triangular and $a_{ij} (D) = 0$ for $i>1$ and $j< n+1$. 

We first compute $B(1+D)B^{-1}$. Since $B = 1 + E_{1,2} + E_{n,n+1}$, we have \begin{align*}
B(1+D)B^{-1} &= (1 + E_{1,2} + E_{n,n+1}) (1+ D)(1 - E_{1,2} - E_{n,n+1})\\
&= 1+  (1 + E_{1,2} + E_{n,n+1}) D(1 - E_{1,2} - E_{n,n+1})\\
&= 1+ D-DE_{1,2} - D E_{n,n+1} + E_{1,2} D + E_{n,n+1} D\\& - E_{1,2} D E_{1,2} -E_{1,2}DE_{n,n+1}-E_{n,n+1}D E_{1,2} -E_{n,n+1}DE_{n,n+1}.
\end{align*}

Since $D$ is strictly upper triangular, $D E_{1,2} = 0$ and $E_{n,n+1} D = 0$. It follows that $$ B(1+D)B^{-1} =  1+ D - D E_{n,n+1} + E_{1,2} D - E_{1,2} D E_{n,n+1}.$$

Since $D_{2,n} = 0$, we have $E_{1,2} D E_{n,n+1} = 0$. Since $D_{j,n} = 0$ for $j>1$, we have that the only non-zero entry of $D E_{n,n+1}$ is $a_{1,n+1} (D E_{n,n+1}) = D_{1,n}$. Since $D_{2,j} = 0$ for $j<n+1$, we have that the only non-zero entry of $E_{1,2} D $ is $a_{1,n+1}(E_{1,2} D) = D_{2,n+1}$. It follows that $$ B(1+D)B^{-1} =  1+ D +(D_{2,n+1}-D_{1,n}) E_{1,n+1}.$$

Therefore, \begin{align*}[B, C] &= (C +(D_{2,n+1}-D_{1,n}) E_{1,n+1}) C^{-1} \\
 & = 1+ ((D_{2,n+1}-D_{1,n}) E_{1,n+1}) C^{-1} = 1+ (D_{2,n+1}-D_{1,n}) E_{1,n+1}\\
 & = 1+ (a_{2,n+1}(C) - a_{1,n} C)E_{1,n+1}, \end{align*} where the second equality follows because $C$ and $C^{-1}$ are upper triangular matrices with all diagonal entries equal to $1$.
\end{proof}

The basis $\{x,y\}$ determines an isomorphism $$\pi/[\pi]_2 \cong \Z^{\Sigma}(\chi)^2,$$ whence an isomorphism $H^1(G, \pi^{ab}) \cong H^1(G, \Z^{\Sigma}(\chi))^2$. For an element $[p]$ of $H^1(G, \pi/[\pi]_{n})$, let $[p]_x \oplus [p]_y $ denote the image of $[p]$ under the map $$ H^1(G, \pi/[\pi]_{n}) \to  H^1(G, \pi^{ab}) \cong H^1(G, \Z^{\Sigma}(\chi))^2.$$ Similarly, for a cocycle $p$ in $C^1(G, \pi/[\pi]_{n})$, let $p_x \oplus p_y $ denote the image of $p$ under the map $$ C^1(G, \pi/[\pi]_{n}) \to  C^1(G, \pi^{ab}) \cong C^1(G, \Z^{\Sigma}(\chi))^2.$$ 

\tpoint{Theorem}\label{Main_thm_obstruction_section}{\em Let $p: G \to \pi/[\pi]_n$ be a cocycle, and let $[p]$ denote the corresponding cohomology class. If $\delta_n[p] = 0$, then the order $n$-Massey product $\langle -[p]_y, - [p]_x, \ldots, - [p]_x , - [p]_y \rangle$ containing two copies of $-[p]_y$ and $(n-1)$ copies of $-[p]_x$  admits the defining system determined by $\varphi \circ p$, and with respect to this defining system $$\langle -[p]_y, - [p]_x, - [p]_x, \ldots, - [p]_x, -[p]_x , -  [p]_y \rangle = - [\h] \cup [p]_y,$$ where $[\h]$ in $H^1(G, \Z^{\Sigma}(\chi^{n-1}))$ is represented by the cocycle $$\h(g) = a_{2,n+1}(\varphi (\mathfrak{f}(g))) - a_{1,n}(\varphi (\mathfrak{f}(g))).$$}

\begin{proof}
Suppose $\delta_n[p] = 0$. Then, there exists a cocycle $q': G \to \pi/[\pi]_{n+1}$ whose composition $q$ of $q'$ with the quotient map $\pi/[\pi]_{n+1} \to \pi/[\pi]_n$ is cohomologous to $p$. This implies that there exists $\gamma$ in $\pi/[\pi]_n$ such that $\gamma^{-1} q(g) (g\gamma) = p(\gamma)$ for all $g$ in $G$ \cite[p. 16]{coh_num_fields}. Choose $\tilde{\gamma}$ in $\pi/[\pi]_{n+1}$ mapping to $\gamma$ under the quotient map. By replacing $q'$ with $g \mapsto \tilde{\gamma}^{-1} q'(g) (g\tilde{\gamma})$, we may assume that $q = p$.

The morphisms $\varphi$ and $\varphi'$ fit into a commutative diagram of set-maps \begin{equation}\label{varphi'varphiCD}
\xymatrix{&& [\pi]_{n}/[\pi]_{n+1} \ar[d] & \Z^{\Sigma}(\chi^n)\ar[d] \\ && \pi/[\pi]_{n+1} \rtimes G \ar[d] \ar@{.>}[r]^{\varphi' \rtimes 1} & U_{n+1} \rtimes G \ar[d] \\ G\ar[rr]^{q \rtimes 1}\ar[rru]^{q' \rtimes 1} && \pi/[\pi]_{n} \rtimes G \ar[r]^{\varphi \rtimes 1}& \overline{U}_{n+1} \rtimes G}
 ,\end{equation} where the solid arrows are also group homomorphisms. 

Then $\varphi \circ p: G \to \overline{U}_{n+1}$ determines $$(\varphi \circ p \rtimes 1)=(\varphi \rtimes 1) \circ (q \rtimes 1): G \to \overline{U}_{n+1} \rtimes G,$$ which determines the defining system \begin{align}\label{defining_system_<y,x,...,x,y>}\{-a_{i,j} \circ (\varphi \rtimes 1) \circ (q \rtimes 1) : i,j\in\{1,2,\ldots, n+1 \}, i<j\text{ and }(i,j) \neq (1,n+1)\}  \nonumber \\
 = \{-a_{i,j} \varphi p: i,j\in\{1,2,\ldots, n+1 \}, i<j\text{ and }(i,j) \neq (1,n+1)\}
\end{align} for the order $n$ Massey product  
$$\langle -a_{1,2} \varphi p , -a_{2,3} \varphi p , \ldots, -a_{n,n+1}\varphi p \rangle.$$ See Definition \ref{def_defining_system_determined_by_theta}.

By definition of $\varphi$, we have that $$a_{1,2} \varphi p= a_{n,n+1}\varphi p= p_y $$ $$a_{i,i+1} \varphi p = p_x \text{ for } i =2,3,\ldots, n-1.$$ Thus \eqref{defining_system_<y,x,...,x,y>} is a defining system for the order $n$ Massey product $$ \langle -[p]_y, -[p]_x, \ldots, -[p]_x, -[p]_y \rangle.$$

We may evaluate this Massey product with respect to the defining system \eqref{defining_system_<y,x,...,x,y>} as follows.

Let $\omega$ denote the element of $H^2(\overline{U}_{n+1} \rtimes G, \Z^{\Sigma}(\chi^n))$ classifying the short exact sequence in the right vertical column of \eqref{varphi'varphiCD}. Then $(\varphi p \rtimes 1)^* \omega$ classifies the short exact sequence \begin{equation}\label{pull_back_by_petc_SES}
1 \to \Z^{\Sigma}(\chi^n) \to G \times_{\overline{U}_{n+1} \rtimes G} (U_{n+1} \rtimes G) \to G \to 1
\end{equation} obtained by pulling back the right vertical column of \eqref{varphi'varphiCD} by $\varphi p \rtimes 1$. By \cite[Remark p. 182]{Dwyer} or see \cite[2.4]{Kyoto} for a modification for profinite groups with actions, $\omega$ is the Massey product $\langle -a_{1,2}, -a_{2,3}, \ldots, -a_{n,n+1} \rangle$ with respect to the defining system $-a_{i,j}$. It follows that the element $(\varphi p \rtimes 1)^* \omega$ is the Massey product $$ \langle -[p]_y, -[p]_x, \ldots, -[p]_x, -[p]_y \rangle$$ with respect to the defining system \eqref{defining_system_<y,x,...,x,y>}. 

By \eqref{varphi'varphiCD}, $s=(\varphi' \rtimes 1) \circ (q' \rtimes 1)$ is a continuous section of the quotient map of the short exact sequence \eqref{pull_back_by_petc_SES}. Therefore, the cocycle $\tilde{s}$ determined by $s$ via the formula \begin{equation*}\tilde{s}(g,h)= s(g)s(g)s(gh)^{-1} \end{equation*} represents $(\varphi p \rtimes 1)^* \omega$ by \cite[IV \S 3 (3.3)]{Brown_coh_groups}.

Let $\tilde{\varphi} =\varphi \circ q'$, so $s(g) = \tilde{\varphi} g \rtimes g$. Thus \begin{align*}
\tilde{s}(g,h) &= (\tilde{\varphi}(g) \rtimes g)(\tilde{\varphi}(h) \rtimes h)(\tilde{\varphi}(gh) \rtimes gh)^{-1}\\
& = \tilde{\varphi}(g)(g\tilde{\varphi}(h))(\tilde{\varphi}(gh))^{-1}. 
\end{align*} 

Since $q'$ is a twisted homomorphism,  $q'(gh) = q'(g) (g q'(h))$. Since $\varphi'$ is a homomorphism $$\tilde{\varphi}(gh) = \varphi' q'(gh) = \varphi' q'(g) \varphi'(g q'(h)).$$ Thus \begin{align*}
\tilde{s}(g,h) =  \tilde{\varphi}(g)(g\varphi' q'(h))(\varphi'(g q'(h)))^{-1} (\tilde{\varphi}(g))^{-1}. 
\end{align*}
Since the quotient map $U_{n+1} \to \overline{U}_{n+1}$ is $G$-equivariant, and so is $\varphi'$ composed with this quotient map, it follows that $(g\varphi' q'(h))(\varphi'(g q'(h)))^{-1}$ is in the kernel of $U_{n+1} \to \overline{U}_{n+1}$. It follows that $(g\varphi' q'(h))(\varphi'(g q'(h)))^{-1}$ is in the center of $U_{n+1}$, and thus commutes with  $\tilde{\varphi}(g)$. Thus $$\tilde{s}(g,h) = (g\varphi' q'(h))(\varphi'(g q'(h)))^{-1}.$$

Because there is a continuous set-theoretic section of $\pi \to \pi/[\pi]_2 \cong  (\Z^{\Sigma})^2$, we may write $q'(h)$ as $$q'(h)= y^{q'_y(h)} x^{q'_x(h)} q'_r(h)$$ with $q'_r(h)$ in $[\pi]_2$, $q'_x(h)$ and $q'_y(h)$ in $\Z^{\Sigma}$. Thus $g q'(h) = gy^{q'_y(h)} gx^{q'_x(h)}  gq'_r(h)$ and $\varphi'(g q'(h)) = \varphi'(gy^{q'_y(h)}) \varphi'(gx^{q'_x(h)}) \varphi'(gq'_r(h))$.

For all $g$ in $G$, $r$ in $[\pi]_2$, and $e \in \Z^{\Sigma}$ we claim that \begin{enumerate}
\item \label{varphigr=gvarphir}$\varphi'(gr) = g \varphi' (r).$
\item \label{varphigx=gvarphix} $\varphi'(g x^e) = g \varphi'(x^e).$
\end{enumerate}

To see \eqref{varphigr=gvarphir}, note that for all $a,b$ in $\pi/[\pi]_{n+1}$, \begin{align*}(g\varphi'(ab))(\varphi'(g(ab)))^{-1} &= g \varphi'(a) g \varphi'(b) (\varphi'(ga) \varphi'(gb))^{-1}  \\
&= g \varphi'(a) g \varphi'(b) ( \varphi'(gb))^{-1} ( \varphi'(ga))^{-1} \\
&= g \varphi'(a) ( \varphi'(ga))^{-1}  g \varphi'(b) ( \varphi'(gb))^{-1}\end{align*} where the second equality follows because $g \varphi'(b) (g \varphi'(b))^{-1}$ is in the center of $U_{n+1}$. In other words $a \mapsto g \varphi'(a) ( \varphi'(ga))^{-1} $ is a homomorphism $\pi/[\pi]_{n+1} \to \Z^{\Sigma}$. Since $\Z^{\Sigma}$ is abelian, this map must factor through $\pi/[\pi]_1$, which shows  \eqref{varphigr=gvarphir}. 

To see \eqref{varphigx=gvarphix}, note that \begin{align*} \varphi'(g x^e) = \varphi'((gx)^e) = \varphi'(gx)^e .\end{align*} By Lemma \ref{gvarphi'varphi'gdif}, $\varphi'(gx) = g \varphi'(x)$, whence \begin{align*} \varphi'(gx)^e = (g \varphi'(x))^e = g \varphi'(x)^e  = g \varphi'(x^e).\end{align*} Combining the two previous, we have $\varphi'(g x^e) =g \varphi'(x^e) $, showing \eqref{varphigx=gvarphix}.

From \eqref{varphigr=gvarphir} and \eqref{varphigx=gvarphix}, it follows that $$\tilde{s}(g,h) =( g \varphi'( y^{q'_y(h)}))(  \varphi'(gy^{q'_y(h)}) )^{-1}.$$

Let $e$ be an element of $\Z^{\Sigma}$. Then, \begin{align*}
\varphi'(g y^e)& = \varphi'(\mathfrak{f}(g)^{-1}y^{e\chi(g)}\mathfrak{f}(g) ) \\
&= \varphi'(\mathfrak{f}(g))^{-1} g \varphi'(y)^e \varphi'(\mathfrak{f}(g)) \\
&=  [\varphi'(\mathfrak{f}(g))^{-1}, g \varphi'(y)^e] g \varphi'(y)^e.
\end{align*} where the second equality follows from the equality $g \varphi'(y)^e = \varphi'(y^{e\chi(g)})$ resulting from Lemma \ref{gvarphi'varphi'gdif}. 

Thus, $$\tilde{s}(g,h) =  [\varphi'(\mathfrak{f}(g))^{-1}, g \varphi'(y)^{q'_y(h)}]^{-1}.$$

By Lemma \ref{varphi(commutator)_in_N}, $\varphi(\gamma)$ is in $\N$ for any $\gamma$ in $[\pi]_2$. Since $\varphi (y)$ is in the center of $\N$, it follows that $\varphi(y)^e$ is in the center of $\N$, as well as $g \varphi(y)^e$ for any $g$. It follows that $[\varphi'(\mathfrak{f}(g))^{-1},\delta]$ is in the kernel of $U_{n+1} \to \overline{U}_{n+1}$, when $\delta = g \varphi'(y)^{q'_y(h)}, g \varphi'(y)$ or $\varphi'(y)$. Since this kernel is in the center, it follows that $[\varphi'(\mathfrak{f}(g))^{-1},\delta]$ is in the center of $U_{n+1}$. By the Witt-Hall identity $[a,bc]= [a,b][a,c][[c,a],b]$\hidden{proof: [a,b][a,c][[c,a],b] = [a,b]b[a,c] b^{-1} = a b a^{-1} a c a^{-1} c^{-1} b^{-1} = a b c a^{-1} c^{-1} b^{-1} = [a,bc]} and the equality $$g \varphi'(y) = \varphi'(y)^{\chi(g)} ,$$ it follows that  \begin{align*}[\varphi'(\mathfrak{f}(g))^{-1}, g \varphi'(y)^{q'_y(h)}]^{-1} &= [\varphi'(\mathfrak{f}(g))^{-1},\varphi'(y)]^{-\chi(g)q'_y(h)} \\
&= [\varphi'(y), \varphi'(\mathfrak{f}(g))^{-1}]^{\chi(g)q'_y(h)} \\
&= [\varphi'(y), \varphi'(\mathfrak{f}(g))]^{-\chi(g)q'_y(h)}.\end{align*}

Thus, \begin{align*}\tilde{s}(g,h) & = [B,\varphi'(\mathfrak{f}(g)) ]^{-\chi(g)q'_y(h)} \\ &= a_{1,n+1}( [B,\varphi'(\mathfrak{f}(g)) ])(-\chi(g)q'_y(h)). \end{align*}

By Lemma \ref{[B,C]CinN}, \begin{align*}a_{1,n+1}( [B,\varphi'(\mathfrak{f}(g)) ])& = a_{2,n+1}(\varphi'(\mathfrak{f}(g)) ) - a_{1,n} \varphi'(\mathfrak{f}(g))\\ 
& = a_{2,n+1}(\varphi(\mathfrak{f}(g)) ) - a_{1,n} \varphi(\mathfrak{f}(g)),\end{align*} proving the theorem.
\end{proof}

\epoint{Remark}\label{main_thm_mod_remark} We could instead define $U_{n+1}$ and $\overline{U}_{n+1}$ to have coefficients in $\Z^{(p')}[\frac{1}{n!}]$, where $\Z^{(p')}$ denotes the prime-to-p completion of $\Z$, i.e., the inverse limit of all $\Z/m$ with $m$ not divisible by $p$. We could then let $\pi$ be the prime-to-p completion of the free group $\langle x,y\rangle$. Theorem \ref{Main_thm_obstruction_section} remains valid when we interpret $[p]_y$, $[p]_x$ as elements of $H^1(G,\Z^{(p')}[\frac{1}{n!}](\chi) )$ and $[\h]$ as an element of $H^1(G, \Z^{(p')}[\frac{1}{n!}](\chi^{n-1}))$.

\epoint{Remark}\label{main_theorem_as_partial_computation_delta_n} Theorem \ref{Main_thm_obstruction_section} can be restated as a partial computation of Ellenberg's obstruction $\delta_n$. The homomorphism $\varphi'$, induces a homomorphism $\mu: [\pi]_n / [\pi]_{n+1} \to \Z^{\Sigma}(\chi^n)$. It follows from \cite[Lem 4.2]{Dwyer} that $\mu(\gamma)$ is the Magnus coefficient $\mu(y,x,x,\ldots,x,x,y; \gamma)$. Magnus coefficients and the Magnus embedding are recalled in \cite[2.5]{Kyoto}. The Magnus embedding determines an isomorphism of $[\pi]_n / [\pi]_{n+1}$ into the homogeneous degree $n$ Lie elements of the non-commutative power series ring in the variables $x$ and $y$ \cite[\S 5.7,~Cor.~5.12(i)]{Magnus_Karrass_Solitar}, and it follows that the direct sum of the degree $n$ Magnus coefficients determines an injective equivariant homomorphism, which is also the inclusion of a direct summand. (This is discussed in detail in \cite[2.5]{Kyoto}.) It follows that the computation of $\mu'_* \delta_n$ for all degree $n$ Magnus coefficients $\mu'$ gives the computation of $\delta_n$. In \cite{Kyoto}, there is a computation of $\mu'_* \delta_n$ when $\mu'$ is a degree $n$ Magnus coefficient with one appearance of the variable $y$. Theorem  \ref{Main_thm_obstruction_section} computes $\mu_* \delta_n$ for the Magnus coefficient $\mu$.  To see this, let $\omega$ be the element of $H^2(\overline{U}_{n+1} \rtimes G, \Z^{\Sigma}(\chi^n))$ classifying the short exact sequence in the right vertical column of \eqref{varphi'varphiCD} as above. Let $\eta$ denote the element of $H^2( \pi/[\pi]_{n} \rtimes G,  [\pi]_{n}/[\pi]_{n+1})$ classifying the short exact sequence in the left vertical column of \eqref{varphi'varphiCD}. It is tautological that $\mu_* \delta_n [p] = p^* \mu_* \eta$. Let $r$ be a continuous section of the quotient map $\pi/[\pi]_{n+1} \to \pi/[\pi]_n$. Then $r \rtimes 1$ is a continuous section of $\pi/[\pi]_{n+1} \rtimes G \to \pi/[\pi]_n \rtimes G$. Using the section $r \rtimes 1$ to compute a representative of $p^*\mu_* \eta$ produces the cocycle $$(g,h) \mapsto \varphi'(r(p(g))) \varphi' (g (r(p(h))) ) \varphi' (r(p(gh)))^{-1}.$$ Using the section $\varphi' r \rtimes 1$ to compute a representative of $$p^* (\varphi \rtimes 1)^* \omega = \langle -[p]_y, - [p]_x, \ldots, - [p]_x , - [p]_y \rangle$$ produces the cocycle $$(g,h) \mapsto \varphi'(r(p(g))) g\varphi' ( (r(p(h))) ) \varphi' (r(p(gh)))^{-1}.$$ Comparing these two cocycles, we compute that $$ \langle -[p]_y, - [p]_x, \ldots, - [p]_x , - [p]_y \rangle - \mu_* \delta_n [p] $$ is represented by the cocycle $$(g,h) \mapsto  g\varphi' ( (r(p(h))) ) \varphi' (g (r(p(h))) ) ^{-1}.$$ The computation in the proof of Theorem \ref{Main_thm_obstruction_section} shows that this cocycle represents $-[\h] \cup [p]_y$. (For this, let $q' = rp$ and note that the necessary part of the computation does not use the fact that $q'$ is a homomorphism.) Thus, the we have \begin{equation}\label{mu_*delta_n_comp} \mu_* \delta_n [p] = [\h] \cup [p]_y +  \langle -[p]_y, - [p]_x, \ldots, - [p]_x , - [p]_y \rangle,\end{equation} where the defining system for the Massey product is as specified in Theorem \ref{Main_thm_obstruction_section}.

\section{Applications}

This section gives two applications of Theorem \ref{Main_thm_obstruction_section}: the first to $\pi_1$-sections of $\mathbb{P}^1 - \{0,1,\infty\}$ and the second to Massey products in Galois cohomology. Except for the statements and proofs of Corollaries \ref{Cor_pi1-sections_Pic_lift_obstruction} and \ref{Cor_MPGA}, the material in this section is expository or establishes notation.

\subsection{Notation}\label{Applications_section_notation_subsection} As in Assumption \ref{assumption_n>=3}, we again fix an integer $n \geq 3$. 

Let $k$ be a field of characteristic $p \geq 0$. Let $\Sigma$ be a set of primes not including $p$ or any prime dividing $n!$. Let $G = \Galk$.

Let $\chi: \Galk \to (\Z^{\Sigma})^*$ denote the cyclotomic character, defined so that the image $\chi(g)_m$ of $\chi(g)$ in $(\Z/m)^*$ is such that $g \zeta_m = \zeta_m^{\chi(g)_m}$.

Let $X$ be a scheme over $k$, equipped with a geometric point. We use this geometric point as a base point for the \'etale fundamental group. Since we could choose a geometric point associated to a rational point or rational tangential base point, and the constructions we then do with this geometric point are independent of the choice, we will also say that $X$ is equipped with a base point, and allow this base point to be a rational point or a rational tangential base point. The kernel $K$ of the map $\pi_1(X_{\kc}) \to \pi_1(X_{\kc})^{\Sigma}$ is a characteristic subgroup. Pushing out the homotopy exact sequence by this kernel yields \begin{equation}\label{piSigma_SES_diagram} \xymatrix{ 1 \ar[r] & \pi_1(X_{\kc}) \ar[d]  \ar[r] & \pi_1(X) \ar[r] \ar[d] & \Galk \ar[d] \ar[r] & 1 \\ 1 \ar[r] & \pi_1(X_{\kc})^{\Sigma} \ar[r] & \pi_1(X)/K \ar[r] & \Galk \ar[r] & 1} \end{equation} A {\em $\pi_1^{\Sigma}$-section} of $X$ means a section of the bottom exact sequence of \eqref{piSigma_SES_diagram}. 

In the case where the bottom exact sequence of \eqref{piSigma_SES_diagram} is already equipped with a fixed splitting, it becomes the sequence \begin{equation*}\label{piSigmasplitSES} 1 \to \pi_1(X_{\kc})^{\Sigma}  \to \pi_1(X_{\kc})^{\Sigma}  \rtimes \Galk \to \Galk  \to 1.\end{equation*} A splitting is then equivalent to a cocycle $$\Galk \to  \pi_1(X_{\kc})^{\Sigma} $$ and by a slight abuse of notation, we will also call such a cocycle a $\pi_1^{\Sigma}$-section. By another slight abuse of notation, the term $\pi_1^{\Sigma}$-section also sometimes refers to the conjugacy class of a splitting or the cohomology class of a cocycle.

\epoint{Kummer map}\label{Kummer_map_recall} Let $r$ be an integer relatively prime to $p$. Let $H^*(\Spec k, -)$ denote \'etale cohomology. Applying $H^*(\Spec k, -)$ to the short exact sequence $$1 \to \mu_r \to \G_{m} \stackrel{z \mapsto z^r}{\to} \G_m \to 1 $$ of group schemes over $k$ gives the exact sequence $$ \ldots \to k^* \stackrel{z \mapsto z^r}{\to} k^* \to H^1(\Galk, \Z/r(\chi)) \to 0 \ldots$$ by Hilbert 90. Allowing $r$ to vary over the integers whose prime factors are in $\Sigma$, we obtain the Kummer map $$k^* \to H^1(\Galk, \Z^{\Sigma}(\chi)) ,$$ and the isomorphism $$(k^*)^{\Sigma} \cong H^1(\Galk, \Z^{\Sigma}(\chi)),$$ where $(k^*)^{\Sigma}$ denotes $\varprojlim_r k^*/(k^*)^r$, where $r$ varies over integers whose prime factors are contained in $\Sigma$.

\subsection{Rational tangential base points}\label{rational_tangential_basepoints_subsection} In \cite[\S 15]{Deligne}, Deligne defines tangential base points for the \'etale fundamental group, primarily for curves over a field of characteristic $0$, but he comments that by restricting to moderately ramified covers, the construction works in characteristic $p>0$ \cite[15.26]{Deligne}. In this subsection, we recall the notion of a tangential base point, introduce some needed notation, and make a few remarks so that we may use tangential base points for $\pi_1^{\et,\Sigma}$ in characteristic $p \geq 0$. We follow the point of view of \cite{Nakamura}, making the necessary modifications for characteristic $p>0$, provided $\Sigma$ does not contain $p$, as we have assumed.

Let $\overline{X}$ be a smooth, proper, geometrically connected curve over $k$, and let $X \subseteq \overline{X}$ be an open subset. Choose $x$ in $\overline{X}(k)$. The completed local ring $\hat{\mathcal{O}}_{\overline{X},x}$ of $\overline{X}$ at $x$ is isomorphic to $k[[z]]$ and a choice of a local parameter $z$ produces such an isomorphism. 

\tpoint{Lemma}{\em The field of Puiseux series $$\kc((z^{\Q})) = \cup_{m =1}^{\infty} \kc((z^{\frac{1}{m}}))$$ contains every finite Galois extension of the Laurent series ring $L((z))$ whose degree is not divisible by $p$ for any extension $L$ of $k$ in $\kc$.}

\begin{proof} cf. \cite[Lemma 3]{Kedlaya_Puiseux}. We may assume that $L=\kc$. If $p=0$, Puiseux's theorem implies $\kc((z^{\Q}))$ is algebraically closed, so we may assume $p>0$. Let $E$ be a Galois extension of degree prime to $p$. The degree of the residue field extension must also be prime to $p$, and therefore the residue field extension is finite and separable, whence trivial. Thus $E$ is totally ramified. Since the wild inertia group is a $p$-group, it is trivial. Since the quotient of the inertia group by the wild inertia group is cyclic of degree $m$ with $m$ prime to $p$, we have that $E$ is a cyclic $m$-extension. Since $m$ is not divisible by $p$, $\kc$ contains the $m$th roots of unity. By Kummer theory, there is $x$ such that $E = \kc((z))(x^{\frac{1}{m}})$. By the equality $(1+b)^{\frac{1}{m}} = \sum_{i=0}^{\infty} {\frac{1}{m} \choose i} b^i$, which is valid in characteristic $p$ for $m$ prime to $p$, it follows that $E$ is contained in $\kc((z^{\frac{1}{m}}))$, completing the proof.
\end{proof}

 Let $\Omega$ be an algebraically closed field extension of $\kc((z^{\Q}))$. The composition $$ \Spec \Omega \to \Spec \kc((z^{\Q})) \to \Spec k[[z]] \cong \hat{\mathcal{O}}_{\overline{X},x} \to \overline{X} $$ factors through the generic point of $\overline{X}$ and therefore determines a geometric point \begin{equation*}\label{bz} b:\Spec \Omega \to \Spec \kc((z^{\Q})) \to X.\end{equation*}  Any such map will be called the {\em tangential base point} of $X$ at $x$ in the direction of $z$, or the {\em rational tangential base point} to emphasize that $x$ was a rational point.  The geometric point $b$ determines an embedding $k(X) \subset k((z)) \subset \kc((z^{\Q}))$. Let $E$ be the fixed field of the separable closure of $\kc((z))$ under the kernel of $\Gal(\kc((z))^s/\kc((z))) \to \Gal(\kc((z))^s/\kc((z)))^{\Sigma}$. Note that $E$ is a subfield of $ \kc((z^{\Q}))$. The coefficientwise action of $\Galk$ on $ \kc((z^{\Q}))$ determines splittings of $\Gal(E/\kc((z^{\Q})) ) \to \Galk$ as well as $\Gal(E(k(X))/k(X)) \to \Galk$, where $E(k(X))$ denotes the fixed field of the separable closure of $\kc(X)$ under the kernel of the map from $\Gal(\kc(X)^s/\kc(X))$ to its pro-$\Sigma$ completion. It then follows \cite[V Prop 8.2]{sga1} that $b$ determines a splitting of the bottom exact sequence of \eqref{piSigma_SES_diagram}, i.e., a   $\pi_1^{\Sigma}$-section of $X$. A local parameter $z$ at $x$ determines a tangent vector $$\Spec k[[z]]/\langle z^2 \rangle \to \overline{X} $$ and the $\pi_1^{\Sigma}$-section associated to the choice of geometric point only depends on this tangent vector. If $x$ is a rational point of $X$, then the $\pi_1^{\Sigma}$-section only depends on the map $$\Spec k[[z]]/\langle z \rangle \to \overline{X} , $$ i.e., the point $x$ itself. A rational tangential base point such that $x$ is in $\overline{X}(k) - X(k)$ will be called a {\em rational tangential point at infinity}.
 
\epoint{Notation for tangential basepoints of $\mathbb P^1 - \{0,1,\infty \}$}\label{notation_for_tangential_base_points}  The scheme $\mathbb{P}^1_k - \{0,1,\infty\} \cong \Spec k[z, \frac{1}{z}, \frac{1}{z-1}]$ is an open subscheme of $\mathbb{A}^1_k \cong \Spec k[z]$. Denote the tangent vector of $\mathbb{A}^1_k$ $$\Spec k[\epsilon]/\langle \epsilon^2 \rangle \to \mathbb{A}^1_k \cong \Spec k[z] $$ $$ a + x \epsilon \leftarrow z$$ by $a + x \epsilon$. A tangential base point associated to this tangent vector will also be denoted by $a + x \epsilon$. For example, this produces rational tangential base points $0 + x \epsilon$ and $1 + x \epsilon$ for $\mathbb P^1 - \{0,1,\infty \}$ for any $x$ in $k^*$. We will also use rational tangential base points associated to tangent vectors at the point $\infty$ of $\mathbb{P}^1 \supset \mathbb P^1 - \{0,1,\infty \}$. For this, let $\iota: \mathbb{P}_k^1 \to \mathbb{P}_k^1 $ denote the map $z \mapsto 1/z$. Then $\iota(0+ x \epsilon)$ are rational tangential base points of $\mathbb P^1 - \{0,1,\infty \}$ for $x$ in $k^*$.
 
\epoint{Description of $\pi$ and its $G$-action}\label{description_pi_and_G_action} There is a lot of very interesting work on the \'etale fundamental group of $\mathbb P^1 - \{0,1,\infty \}$ due to contributions of Anderson, Coleman, Deligne, Ihara, Kaneko, and Yukinari. See for example, \cite[6.3 Thm~p.115]{Ihara_Braids_Gal_grps}. We recall what we use in this subsection. Let $0 + \epsilon$ be the chosen base point for $\mathbb P^1 - \{0,1,\infty \}$, and let $\pi$ be the pro-$\Sigma$ completion of $\pi_1^{\et}(\mathbb P_{\kc}^1 - \{0,1,\infty \})$. By \cite[XIII Corollaire 2.12]{sga1}, $\pi$ is the pro-$\Sigma$ completion on the free group $\langle x, y \rangle$ on two generators $x$ and $y$, where $x$ generates the inertia group of $0$ and $y$ generates the inertia group of $1$. The proof of  \cite[XIII Corollaire 2.12]{sga1} lifts to characteristic $0$ and then compares to $\mathbb{C}$ using specialization morphisms. Since specialization morphisms preserve path composition, we may moreover assume that $y = \wp^{-1} y' \wp$, where $\wp$ is a path from $0 + \epsilon$ to $1 - \epsilon$, and $y'$ generates the inertia group of $\pi_1(\mathbb P_{\kc}^1 - \{0,1,\infty \}, 1 - \epsilon)$ at $1$. Since the image of the inertia group at $0$ under the map $$\pi_1^{\et}(\mathbb P_{\kc}^1 - \{0,1,\infty \}, 0 + \epsilon) \to \pi_1^{\et}(\mathbb P_{\kc}^1 - \{0,1,\infty \}, 1 - \epsilon)$$ induced by the map $z \mapsto 1-z$ is the inertia group at $1$, we may moreover assume that $y'$ is the image of $x$ under this map. It then follows by the argument in \cite[3.9]{Kyoto} that the $G$-action on $\pi$ is of the form $$g(x) =x^{\chi(g)} $$ $$g(y) =\mathfrak{f}(g)^{-1}y^{\chi(g)}\mathfrak{f}(g) $$ where $\mathfrak{f}: G \rightarrow [\pi]_2$ is a cocycle taking values in the commutator subgroup of $\pi$. 

\subsection{$\pi_1$-sections of $\mathbb P^1 - \{0,1,\infty \}$}\label{subsection_pi1_sectionsP1minus}

The generalized Jacobian $\Pic (\mathbb P_k^1 - \{0,1,\infty \})$ is defined to be the connected component containing the identity of the Picard scheme of the initial one-point compactification of $\mathbb P_k^1 - \{0,1,\infty \}$. This compactification of $\mathbb P_k^1 - \{0,1,\infty \}$ is the scheme $\mathbb P_k^1 / 0 \sim 1 \sim \infty$ formed by glueing the points $0$, $1$, and $\infty$ of $\mathbb P_k^1$ to each other. Since $\Jac \mathbb P_k^1$ is a point, a line bundle on $\mathbb P_k^1 / 0 \sim 1 \sim \infty$ is determined by how the fibers of the pull-back to $\mathbb{P}^1_k$ over $0$, $1$, and $\infty$ are glued together. We obtain a non-canonical isomorphism $$\Jac(\mathbb P_k^1 - \{0,1,\infty \}) \cong \G_{m,k} \times \G_{m,k}.$$ For more information on generalized Jacobians, see \cite{alg_grps_class_fields}.

The Abel-Jacobi map based at $0 + \epsilon$ is the map $$\mathbb P_k^1 - \{0,1,\infty \} \to \Jac(\mathbb P_k^1 - \{0,1,\infty \}) \cong \G_{m,k} \times \G_{m,k}$$ given by $$z \mapsto (z, 1-z).$$ The covers $z \mapsto z^n$ produce an isomorphism $\pi_1^{\et, \Sigma} (\G_{m,\kc},1) \cong \Z^{\Sigma}(\chi)$ and the two projections induce an isomorphism \begin{equation}\label{pi1GmtimesGm} \pi_1^{\et, \Sigma} (\G_{m,\kc} \times \G_{m,\kc},1\times 1) \cong \Z^{\Sigma}(\chi) \oplus \Z^{\Sigma}(\chi)\end{equation} \cite[Proposition 4.7]{Orgogozo}. Since these fundamental groups are abelian, there are canonical isomorphisms between these fundamental groups and those based at another base point.  Applying the functor $\pi_1^{\et, \Sigma}( (-) \otimes \kc)$ with the chosen base point and its image produces the map $$\pi \to \Z^{\Sigma}(\chi) \oplus \Z^{\Sigma}(\chi),$$ where the images of $x$ and $y$ form a basis for $\Z^{\Sigma}(\chi) \oplus \Z^{\Sigma}(\chi)$ as a $\Z^{\Sigma}$ module. (To see this, compose the Abel-Jacobi map with the two projections and note that the inclusion $\mathbb P_{\kc}^1 - \{0,1,\infty \} \to \G_{m,\kc}$ preserves the inertia group at $0$.) 

By definition, the isomorphism \eqref{pi1GmtimesGm}, and the tangential base point of $\G_{m,k} \times \G_{m,k}$ given as the image of $0 + \epsilon$, a $\pi_1^{\et, \Sigma}$-section of $\Pic (\mathbb P_k^1 - \{0,1,\infty \})$ is a splitting of the sequence $$1 \to  \Z^{\Sigma}(\chi) \oplus \Z^{\Sigma}(\chi) \to (\Z^{\Sigma}(\chi) \oplus \Z^{\Sigma}(\chi)) \rtimes G \to G \to 1 .$$ The conjugacy classes of such splittings are in bijection with $$H^1(G, \Z^{\Sigma}(\chi) \oplus \Z^{\Sigma}(\chi) )$$ \cite[IV 2 Prop 2.3]{Brown_coh_groups}. Therefore the conjugacy classes of $\pi_1^{\et, \Sigma}$-sections of $\Pic (\mathbb P_k^1 - \{0,1,\infty \})$ can we written $[p] = [p]_x \oplus [p]_y$, where $[p]_x$ and $[p]_y$ are elements of $H^1(G, \Z^{\Sigma}(\chi)) \cong (k^*)^{\Sigma}$, where $(k^*)^{\Sigma}$ is as defined in \ref{Kummer_map_recall}.

\tpoint{Corollary}\label{Cor_pi1-sections_Pic_lift_obstruction}{\em A $\pi_1^{\et,\Sigma}$-section $p$ of $\Pic (\mathbb P_k^1 - \{0,1,\infty \})$ which lifts to a $\pi_1^{\et,\Sigma}$-section of $\mathbb P_k^1 - \{0,1,\infty \}$ satisfies the condition that the Massey product $$\langle-[p]_y, -[p]_x, \ldots, -[p]_x, -[p]_y  \rangle$$ contains  $$0 \text { and }-[\h] \cup [p]_y.$$ Furthermore, if $\tilde{p}$ denotes such a lift, then $\varphi \tilde{p}$ determines a defining system with respect to which $$\langle-[p]_y, -[p]_x, \ldots, -[p]_x, -[p]_y  \rangle = -[\h] \cup [p]_y.$$ }

\begin{proof}
Let $\alpha: \mathbb P_k^1 - \{0,1,\infty \} \to \Pic (\mathbb P_k^1 - \{0,1,\infty \})$ be the Abel-Jacobi map based at $0+ \epsilon$ as described above. The map $\pi_1^{\et, \Sigma}(\alpha, 0 + 1 \epsilon): \pi \rtimes G \to \pi/[\pi]_2 \rtimes G$ is the semi-deirect product of the abelianization $\pi \to \pi/[\pi]_2$ with the identity map on $G$. Let $\tilde{p}: G \to \pi \rtimes G$ denote a lift of the section $p:G \to \pi/[\pi]_2 \rtimes G$ through the quotient map $\pi_1^{\et, \Sigma}(\alpha, 0 + 1 \epsilon)$. Let $\tilde{p}_n: G \to \pi/[\pi]_n \rtimes G $ denote the composition of $\tilde{p}$ with the appropriate quotient, and by a slight abuse of notation, we will also let $\tilde{p}_n$ denote the corresponding cocycle $\tilde{p}_n: G \to \pi/[\pi]_n$. By construction, $\delta_n ([\tilde{p}_n]) = 0$. By  Theorem \ref{Main_thm_obstruction_section}, we have that $\varphi \circ \tilde{p}_n$ is a defining system for $\langle-[p]_y, -[p]_x, \ldots, -[p]_x, -[p]_y  \rangle$ and that the corresponding Massey product is  $$\langle-[p]_y, -[p]_x, \ldots, -[p]_x, -[p]_y \rangle = -[\h] \cup [p]_y.$$

It remains to show that the Massey product $$\langle-[p]_y, -[p]_x, \ldots, -[p]_x, -[p]_y \rangle$$ contains $0$. Because $\langle-[p]_y, -[p]_x, \ldots, -[p]_x, -[p]_y,  \rangle$ contains $-[\h] \cup [p]_y$, it follows that $D([\h] \cup [p]_y) = 0$. Since $D[p]_y = 0$, it follows that $D[\h] \cup [p]_y = 0$, and therefore that $D[\h] = 0$. Let $Z_{i,j}$ denote the defining system determined by $\varphi \circ \tilde{p}_n$, as in Definition \ref{def_defining_system_determined_by_theta}. We saw that $Z_{1,n+1} = p_y$ in the proof of Theorem \ref{Main_thm_obstruction_section}. Then \begin{equation}\label{defining_system_modified_to_get_0}
Z_{1,n}+[\h], Z_{i,j} \text{ for } i,j \in \{1,2,\ldots, n+1 \},~~ i<j \text{ and } (i,j) \neq (1,2)\text{ or }(1,n+1) 
\end{equation} is also a defining system for $\langle-[p]_y, -[p]_x, \ldots, -[p]_x, -[p]_y \rangle$. With respect to this defining system, \begin{align*}\langle-[p]_y, -[p]_x, \ldots, -[p]_x, -[p]_y \rangle &= -[\h] \cup [p]_y + [\h] \cup Z_{n,n+1} \\
&= -[\h] \cup [p]_y + [\h] \cup [p]_y = 0, \end{align*} as desired.
\end{proof}

Although Corollary \ref{Cor_pi1-sections_Pic_lift_obstruction} gives an obstruction to lifting a $\pi_1^{\et,\Sigma}$-section of $\Pic (\mathbb P_k^1 - \{0,1,\infty \})$ to $\mathbb P_k^1 - \{0,1,\infty \}$, it is not clear that it succeeds in obstructing the lifing of any $\pi_1^{\et,\Sigma}$-section of $\Pic (\mathbb P_k^1 - \{0,1,\infty \})$. In conjunction with a similar obstruction, it does succeed in obstructing the lifting of certain $\pi_1^{\et,\Sigma}$-sections for $n=3$. See \cite[Corollary 38]{PIA}. For any $n$, it can be combined with the obstructions described in \cite{Kyoto}. Its purpose here is to contribute to the program of \cite{PIA} \cite{Kyoto} to describe the  $\pi_1^{\et,\Sigma}$-sections of $\mathbb P_k^1 - \{0,1,\infty \}$.

\subsection{Order $n$-Massey products $\langle y, x,x,\ldots,x,x,y \rangle$ in Galois cohomology} We use the rational points and rational tangential points at infinity of $\mathbb P_k^1 - \{0,1,\infty \}$ to show that certain Massey products of the form $\langle y, x,x,\ldots,x,x,y \rangle$ contain $0$. 

In the following corollary, we identify elements of $k^*$ with their images in $H^1(G_k, \Z^{\Sigma}(\chi))$ under the Kummer map. Note that the multiplicative inverse of an element in $k^*$ corresponds to the additive inverse in $H^1(G_k, \Z^{\Sigma}(\chi))$.

\tpoint{Corollary}\label{Cor_MPGA}{\em Let $(x,y)$ in $H^1(G_k, \Z^{\Sigma}(\chi))^2$ be $(x^{-1}, (1-x)^{-1})$ for $x$ in $k^* -\{1\}$, or $(x^{-1},-x^{-1})$ for $x$ in $k^*$. Then the $n$th order Massey product $$\langle y, x, \ldots, x, y \rangle $$ contains $$0 \text { and }-[\h] \cup [p]_y.$$ Here, the Massey product has $(n-2)$ factors of $x$ and two factors of $y$. }

\begin{proof}
By Corollary \ref{Cor_pi1-sections_Pic_lift_obstruction}, it suffices to show that there are $\pi_1^{\Sigma}$-sections $[\tilde{p}]$ and $[\tilde{q}]$ of $\mathbb P_k^1 - \{0,1,\infty \}$ such that the associated sections $[p]$ and $[q]$ of $\Pic (\mathbb P_k^1 - \{0,1,\infty \})$ satisfy \begin{enumerate}
\item \label{p_rational_pt_x} $[p]_x = x$ and $[p]_y = 1-x$
\item \label{q_rational_tgt_pt_infty} $[q]_x = x$ and $[q]_y = -x$.
\end{enumerate} By the argument of \cite[12.2.2]{PIA} or \cite[Lemma 8]{Kyoto}, the rational point $x$ in $(\mathbb P_k^1 - \{0,1,\infty \} )(k) = k^* -\{1\}$ gives rise to a section $[p]$ satisfying \eqref{p_rational_pt_x}. By the same argument, for $x$ in $k^*$, the rational tangential point $\iota(0 + \frac{1}{x} \epsilon)$ gives rise to a section $[q]$ satisfying \eqref{q_rational_tgt_pt_infty}.
\end{proof}

\epoint{Remark} The defining systems implicit in the conclusion of Corollary \ref{Cor_MPGA} can be made entirely explicit: By Corollary \ref{Cor_pi1-sections_Pic_lift_obstruction}, the Massey product $\langle y, x, \ldots, x, y \rangle$ for $(x,y) = (x^{-1}, (1-x)^{-1})$ with respect to the defining system determined by $\varphi \tilde{p}$ is $$\langle y, x, \ldots, x, y \rangle = -[\h] \cup y $$ and with respect to the defining system \eqref{defining_system_modified_to_get_0}, the same Massey product is $0$. Replacing $p$ by $q$ gives the defining systems for $(x,y) = (x^{-1},-x^{-1})$. 

\epoint{Remark}\label{corMPGA_mod_remark} The conclusion of Corollary \ref{Cor_MPGA} remains valid when elements of $k^*$ are identified with their images in $H^1(\Z^{(p')}[\frac{1}{n!}])$.

%%%%%%%%%%%%%%%%%%%%%%%%%%%%%%%%%
% References
%%%%%%%%%%%%%%%%%%%%%%%%%%%%%%%%%

\bibliographystyle{MPGC}

\bibliography{MPGC}

\end{document}